\documentclass[10pt]{article}

\usepackage{graphics}
\usepackage{graphicx}

\usepackage[a4paper, left=35mm,right=35mm,top=34mm,bottom=34mm]{geometry}
\usepackage[utf8]{inputenc}
\usepackage[T1]{fontenc}
\usepackage[english]{babel}

\usepackage{enumerate}
\usepackage{graphicx}
\usepackage{hyperref}
\hypersetup{
    colorlinks=true,
    linkcolor=blue,
    filecolor=magenta,      
    urlcolor=cyan,
}
\usepackage{listings}
\usepackage{color}

\definecolor{dkgreen}{rgb}{0,0.6,0}
\definecolor{gray}{rgb}{0.5,0.5,0.5}
\definecolor{mauve}{rgb}{0.58,0,0.82}
\lstdefinelanguage{MRGC++}{%
  language=C++,
  morekeywords={T, U, MPI_Irecv, MPI_Isend, MPI_Allreduce, MPI_Waitall, Compute, Map, abs, max, Swap, MPI_Recv_init, MPI_Send_init, MPI_Startall, Copy, Init, InitRecv, InitSend, InitAllReduce, Send, Recv, AllReduce, Finalize, InitSnapshot, Snapshot, SwitchAsync, SnapReduce, MPI_Test, MPI_Start}
}
\usepackage{mathtools,amsthm,amssymb,amsfonts}
\usepackage{algorithm}
\usepackage{algorithmic}
\makeatother
\theoremstyle{plain}

\theoremstyle{definition}

\theoremstyle{remark}

\usepackage{caption} 
\usepackage{fancyhdr}

\author{
  {\normalsize Guillaume Gbikpi-Benissan}\thanks{Universit\'e Paris-Saclay, CentraleSup\'elec, France
    (guibenissan@gmail.com).}
  \and
  {\normalsize Fr\'ed\'eric Magoul\`es}\thanks{Universit\'e Paris-Saclay, CentraleSup\'elec, France (frederic.magoules@hotmail.com).}
}
\title{Accurate Coarse Residual for Two-Level Asynchronous Domain Decomposition Methods}
\date{}

\begin{document}
\maketitle
\thispagestyle{fancy}

\begin{abstract}
\noindent Recently, asynchronous coarse-space correction has been achieved within both the overlapping Schwarz and the primal Schur frameworks. Both additive and multiplicative corrections have been discussed. In this paper, we address some implementation drawbacks of the proposed additive correction scheme. In the existing approach, each coarse solution is applied only once, leaving most of the iterations of the solver without coarse-space information while building the right-hand side of the coarse problem. Moreover, one-sided routines of the Message Passing Interface (MPI) standard were considered, which introduced the need for a sleep statement in the iterations loop of the coarse solver. This implies a tuning of the sleep period, which is a non-discrete quantity. In this paper, we improve the accuracy of the coarse right-hand side, which allowed for more frequent corrections. In addition, we highlight a two-sided implementation which better suits the asynchronous coarse-space correction scheme. Numerical experiments show a significant performance gain with such increased incorporation of the coarse space.
\end{abstract}


\section{Introduction}
\label{sec:intro}

As the amount of parallelism in modern computing platforms is tremendously increasing, numerical iterative methods have to cope with resource failures, heterogeneity and costly communication. Domain decomposition methods \cite{TosWid2005, DoleanEtAl2015} are so far among the approaches providing highest levels of parallelism. Still, their classical mathematical design relies on iterative schemes which hardly support failures and heterogeneity.
Asynchronous iterations \cite{ChazMir1969} early arose as a parallel mathematical scheme where computation and communication can occur concurrently, and are now gaining a particular attention as a viable option for overcoming the limits of classical iterative methods.
For an overview of the asynchronous iterations theory and applications, we refer to, e.g., the review papers \cite{FromSzyld2000, Spiteri2020}.
Recent results have shown promising performances against Krylov methods within a framework of non-overlapping domain decomposition (see, e.g., \cite{GBenissanEtAl2022b, GBenMag2023}). These results were obtained while applying asynchronous coarse-space correction \cite{GlusaEtAl2020b, GBenMag2022}, a recent breakthrough in the asynchronous computing field.

In \cite{GlusaEtAl2020b}, an asynchronous additive coarse-space correction was proposed in the framework of the restricted additive Schwarz (RAS) method, however, an in-depth analysis by \cite{GBenMag2022} in case of multiplicative correction showed that more attention should be paid to the accuracy of the computed coarse residual, which determines the effectiveness of the correction.
A connection can therefore be made between the asynchronous coarse-space correction issue and the asynchronous convergence detection one, both linked by the fact that it seems not trivial to compute a consistent global residual without pausing the ongoing computation. The convergence detection issue has been largely discussed in the literature of asynchronous methods (see, e.g., \cite{SavBert1996, BahiEtAl2008, MielEtAl2008, MagGBen2018b, GBenMag2020}), however, most of the proposed solutions rely on local residuals which, in the context of asynchronous iterations, are not local components of a global residual. In \cite{GlusaEtAl2020b, GlusaEtAl2020}, an approximation approach was used in a centralized form, which only needs a reduction operation over the local residuals (also see, e.g., Fig. 3 or Fig. 4 in \cite{GBenMag2020} or Algorithm 4.5 in \cite{GBenMag2022}
for the distributed form).
A similar operation was performed to also construct the coarse global residual vector. While such a minimal approach can be robust enough for convergence detection, it seems natural that the accuracy of a global residual is much more critical for effective asynchronous coarse-space correction. In this paper, we apply more accurate global residual ideas from \cite{MagGBen2018b, MagGBen2018c} (also see Section 3 of \cite{SavBert1996}) to improve the quality of the coarse-space correction.

The paper is organized as follows. Section~\ref{sec:bg} recalls the iterative model of an RAS solver, its asynchronous implementation, and their current two-level counterparts with additive coarse-space correction. Our contribution is in Section~\ref{sec:new}. We first present a simpler way to manage communications between the subdomains and the coarse domain, then we propose the use of a globally consistent coarse residual to improve the accuracy of the coarse solution. Experimental investigation is conducted in Section~\ref{sec:exp}, and our conclusions follow in Section~\ref{sec:conclusion}.

\section{Computational background}
\label{sec:bg}

We consider iterative solution of a linear algebraic problem
$
A x = b
$
with $A \in \mathbb{R}^{n \times n}$,
where $A$ is a large sparse real coefficients matrix and $x$ is a real unknowns vector.
Consider $p \in \mathbb{N}$ discrete subdomains, and let $R_{i}$, $i \in \{1, \ldots, p\}$, denote the unknowns selection mapping to form the $i$-th subdomain. Derive the local matrices and vectors
$A_{i}^{} = R_{i}^{} A R_{i}^{\mathsf T}$, $x_{i}^{} = R_{i}^{} x$ and $b_{i}^{} = R_{i}^{} b$,
where the superscript $\mathsf T$ denotes the transpose of a matrix or vector. Consider diagonal Boolean matrices, say $B_{i}$, $i \in \{1, \ldots, p\}$, such that
$\sum_{i=1}^{p} R_{i}^{\mathsf T} B_{i} R_{i} = I$,
where $I$ denotes the identity matrix. The RAS preconditioner~\cite{CaiSarkis1999} is then given by
$
M = \sum_{i=1}^{p} R_{i}^{\mathsf T} B_{i}^{} A_{i}^{-1} R_{i}^{},
$
which leads to the RAS iterative solver,
\begin{equation}
\label{eq:relaxation}
x_{}^{k+1} = x_{}^{k} + M \left(b - A x_{}^{k}\right), \qquad k \in \mathbb{N},
\end{equation}
and its parallel form,
\[
x_{i}^{k+1} = \sum_{j=1}^{p} R_{i}^{} R_{j}^{\mathsf T} B_{j}^{} \left(x_{j}^{k} + A_{j}^{-1} \left(b_{j}^{} - R_{j}^{} A x_{}^{k}\right)\right).
\]
Note that $R_{j}^{} A$, $j \in \{1, \ldots, p\}$, is of the form
$
R_{j}^{} A = (A_{j}^{}, R_{j}^{} A R_{\neg j}^{\mathsf T}),
$
where $\neg j$ denotes the complement of the $j$-th subdomain, which means that $R_{j}^{} A R_{\neg j}^{\mathsf T}$ corresponds to the coefficients connecting the unknowns of the $j$-th subdomain to the ones not included in the subdomain. They can be thought of as boundary unknowns. It yields the final parallel scheme
\begin{equation}
\label{eq:ras}
x_{i}^{k+1} = \displaystyle\sum_{j=1}^{p} R_{i}^{} R_{j}^{\mathsf T} B_{j}^{} A_{j}^{-1} \left(b_{j}^{} - R_{j}^{} A R_{\neg j}^{\mathsf T} x_{\neg j}^{k}\right), \quad i \in \{1, \ldots, p\}, \quad k \in \mathbb{N}.
\end{equation}
With such a parallel scheme, communication is needed for boundary components $x_{\neg i}^{k}$ only.

Asynchronous iterations \cite{ChazMir1969} consist of not waiting for sending and receiving data before resuming to the next iteration.
Theoretical modeling and convergence analysis of asynchronous iterations have been largely discussed in the literature. We refer to, e.g., the review papers \cite{FromSzyld2000, Spiteri2020} and the RAS-related papers \cite{FromEtAl1997, FrommerSzyld2001}.
Let $N_{i}$ denote the subset of the subdomains with which the $i$-th subdomain shares boundary unknowns.
Using the well-known Message Passing Interface (MPI) formalism, an asynchronous implementation of the RAS iterations \eqref{eq:ras} is described in Algorithm~\ref{algo:ras:async}.
\begin{algorithm}[htbp]
\caption{Asynchronous RAS solver.}
\label{algo:ras:async}
{\footnotesize
\begin{algorithmic}[1]
\REQUIRE $A_{i}^{}$, $R_{i}^{} A R_{\neg i}^{\mathsf T}$, $b^{}_{i}$, $x_{i}^{}$, $x_{\neg i}^{}$, $B_{i}^{}$, $N_{i}$, $\varepsilon$, $\|b\|$, $k_{\text{max}}$
\STATE $r_{i}^{}$ := $b_{i}^{} - A_{i}^{} x_{i}^{} - R_{i}^{} A R_{\neg i}^{\mathsf T} x_{\neg i}^{}$
; \ AllReduce($r_{i}^{\mathsf T} B_{i}^{} r_{i}^{}$, $r^{\mathsf T} r$, SUM)
; \ $\|r\|$ := $\displaystyle\sqrt{r^{\mathsf T} r}$
; \ $k$ := 0
\STATE req$_{r}$ := REQUEST\_NULL
\WHILE{$\|r\| \ge \varepsilon \|b\|$ \AND $k < k_{\text{max}}$}
	\STATE $x^{}_{i}$ := $x^{}_{i}$ + Solve($A_{i}^{}$, $r^{}_{i}$)
	; \ reqs$_{x}$[ ] := ISynchronize($x_{i}^{}$, $x_{\neg i}^{}$, $N_{i}$)
	; \ FreeAll(reqs$_{x}$)
	\STATE $r_{i}^{}$ := $b_{i}^{} - A_{i}^{} x_{i}^{} - R_{i}^{} A R_{\neg i}^{\mathsf T} x_{\neg i}^{}$
	\IF{Test(req$_{r}$)}
		\STATE $\|r\|$ := $\displaystyle\sqrt{r^{\mathsf T} r}$
		; \ req$_{r}$ := IAllReduce($r_{i}^{\mathsf T} B_{i}^{} r_{i}^{}$, $r^{\mathsf T} r$, SUM)
		; \ $k$ := $k + 1$
	\ENDIF
\ENDWHILE
\RETURN $x^{}_{i}$
\end{algorithmic}
}
\end{algorithm}
The \textit{ISynchronize} function implements necessary communication between the $i$-th subdomain and its neighboring subdomains, in a non-blocking way.
For simplicity, we consider marking the returned requests for being automatically freed upon completion (function \textit{FreeAll}).
Efficient management of MPI communication requests was discussed in, e.g., \cite{ChauEtAl2007, MagGBen2018c}.
In case of an underlying network protocol supporting remote direct memory access (RDMA), one-sided MPI communication should be considered (see, e.g., \cite{WPouChow2016}).

We rather focus on the simple approach analyzed in \cite{GBenMag2020}, where convergence can be reliably detected by only considering the non-blocking collective function \textit{IAllReduce}. It computes an approximation of the dot product $r^{\mathsf T} r$, based on inconsistent local contributions $r_{i}^{\mathsf T} B_{i}^{} r_{i}^{}$. Upon completion of the collective operation (checked with \textit{Test(req$_{r}$)}), all the processes consider a same approximation of the residual norm $\|r\| = \sqrt{r^{\mathsf T} r}$.
While such an approximated residual can be satisfactory in practice for stopping an asynchronous iterative solver, we shall see that it prevents the solver from taking full advantage of coarse-space correction.

In \cite{GlusaEtAl2020b}, an asynchronous implementation of a two-level RAS solver was proposed.
The considered classical coarse-space technique consists of defining a restriction mapping, say $R_{0}^{}$, which aggregates the unknowns of each subdomain into one coarse unknown, and is used to define a coarse matrix, say 
$
A_{0}^{} = R_{0}^{} A R_{0}^{\mathsf T}.
$
The iterative method~\eqref{eq:relaxation} is then modified such as
\[
x_{}^{k+1} = x_{}^{k} + \frac{1}{2} \left(M + R_{0}^{\mathsf T} A_{0}^{-1} R_{0}^{}\right) \left(b - A x_{}^{k}\right), \qquad k \in \mathbb{N},
\]
which results in the two-level parallel scheme
\begin{equation}
\label{eq:ras:cs}
\left\{
\begin{array}{lcl}
x_{0}^{k} & = & A_{0}^{-1} \displaystyle\sum_{j=1}^{p} R_{0}^{} R_{j}^{\mathsf T} B_{j}^{} \left(b_{j} - A_{j}^{} x_{j}^{k} - R_{j}^{} A R_{\neg j}^{\mathsf T} x_{\neg j}^{k}\right),\\
x_{i}^{k+1} & = & \displaystyle\sum_{j=1}^{p} \dfrac{1}{2} R_{i}^{} R_{j}^{\mathsf T} B_{j}^{} \left(A_{j}^{-1} \left(b_{j}^{} - R_{j}^{} A R_{\neg j}^{\mathsf T} x_{\neg j}^{k}\right) + x_{j}^{k} + R_{j}^{} R_{0}^{\mathsf T} x_{0}^{k}\right),
\end{array}
\right.
\end{equation}
with $i \in \{1, \ldots, p\}$ and $k \in \mathbb{N}$. A comprehensive introduction to multilevel techniques can be found, e.g., in the review books \cite{Hackbusch1985, TosWid2005, DoleanEtAl2015}.

The asynchronous implementation proposed in \cite{GlusaEtAl2020b} is based on the fact that the coarse domain is simply seen as a $(p+1)$-th subdomain handling the component $x_{0}^{}$. The asynchronous implementation therefore consists of not waiting for an updated $x_{0}^{}$ before updating $x_{i}^{}$. However, the coarse solver waits for receiving each of the contributions $R_{0}^{} R_{i}^{\mathsf T} B_{i}^{} r_{i}^{}$ before computing $x_{0}^{}$. Similarly, the update of $x_{i}^{}$ is corrected with $x_{0}^{}$ only when an updated $x_{0}^{}$ is received.
It follows that a coarse solution $x_{0}^{}$ can be rarely available, and when it is finally provided, it is applied only once. We propose here a slight modification of this implementation approach, which allowed us to more often correct $x_{i}^{}$, whether $x_{0}^{}$ has been updated or not.

\section{New computational model}
\label{sec:new}

\subsection{Two-level asynchronous RAS solver with two-sided MPI}

The two-level asynchronous RAS solver proposed in \cite{GlusaEtAl2020b, GlusaEtAl2020} is based on MPI one-sided communication routines, however, coarse-space correction requires checking arrival of new data. As a drawback, the coarse domain process has to constantly poll shared Boolean variables in an active loop, waiting for subdomain contributions to the coarse residual vector. A sleep statement was therefore introduced in the loop, hence, one has to determine the optimal sleep period. This adds to the other problem of using each coarse solution only once.

Two-sided routines are particularly suitable for such data arrival checks. If a $(p+1)$-th process is dedicated to the coarse domain, a function \textit{Wait} can be simply called on requests initialized by \textit{IRecv} or \textit{IGather}. If, otherwise, a $i_{0}$-th process, $i_{0} \in \{1, \ldots, p\}$, handles the coarse domain, then completion is just checked with a function \textit{Test}. The latter case is described here in Algorithm~\ref{algo:ras:async:cs}.
\begin{algorithm}[htbp]
\caption{Two-level asynchronous RAS with non-blocking collective routines.}
\label{algo:ras:async:cs}
{\footnotesize
\begin{algorithmic}[1]
\REQUIRE $A_{i}^{}$, $R_{i}^{} A R_{\neg i}^{\mathsf T}$, $b^{}_{i}$, $x_{i}^{}$, $x_{\neg i}^{}$, $B_{i}^{}$, $N_{i}$, $A_{0}^{}$, $R_{0}^{} R_{i}^{\mathsf T}$, $i_{0}^{}$, $\varepsilon$, $\|b\|$, $k_{\text{max}}$
\STATE $r_{i}^{}$ := $b_{i}^{} - A_{i}^{} x_{i}^{} - R_{i}^{} A R_{\neg i}^{\mathsf T} x_{\neg i}^{}$
; \ AllReduce($r_{i}^{\mathsf T} B_{i}^{} r_{i}^{}$, $r^{\mathsf T} r$, SUM)
; \ $\|r\|$ := $\displaystyle\sqrt{r^{\mathsf T} r}$
; \ $k$ := 0
\STATE req$_{r}$ := REQUEST\_NULL
; \ state$_{0}$ := 0
\WHILE{$\|r\| \ge \varepsilon \|b\|$ \AND $k < k_{\text{max}}$}
	\STATE \textbf{if} state$_{0}$ = 0 \textbf{then}
		req$_{r_{0}^{}}$ := IReduce($R_{0}^{} R_{i}^{\mathsf T} B_{i}^{} r_{i}^{}$, $r_{0}^{}$, $i_{0}^{}$, SUM)
		; \ state$_{0}$ := 1
	\textbf{end if}
	\IF{state$_{0}$ = 1 \AND Test(req$_{r_{0}^{}}$)}
		\STATE \textbf{if} $i$ = $i_{0}$ \textbf{then}
			$x_{0}$ := Solve($A_{0}^{}$, $r^{}_{0}$)
		\textbf{end if}
		\STATE req$_{x_{0}^{}}$ := IBcast($x_{0}^{}$, $i_{0}^{}$)
		; \ state$_{0}$ := 2
	\ENDIF
	\STATE \textbf{if} state$_{0}$ = 2 \AND Test(req$_{x_{0}^{}}$) \textbf{then}
		$x^{}_{i}$ := $x^{}_{i}$ + 0.5 $\times$ (Solve($A_{i}^{}$, $r^{}_{i}$) + $R_{i}^{} R_{0}^{\mathsf T} x_{0}^{}$)
		; \ state$_{0}$ := 0
	\STATE \textbf{else}
		$x^{}_{i}$ := $x^{}_{i}$ + Solve($A_{i}^{}$, $r^{}_{i}$)
	\textbf{end if}
	\STATE reqs$_{x}$[ ] := ISynchronize($x_{i}^{}$, $x_{\neg i}^{}$, $N_{i}$)
	; \ FreeAll(reqs$_{x}$)
	; \ $r_{i}^{}$ := $b_{i}^{} - A_{i}^{} x_{i}^{} - R_{i}^{} A R_{\neg i}^{\mathsf T} x_{\neg i}^{}$
	\STATE \textbf{if} Test(req$_{r}$) \textbf{then}
		$\|r\|$ := $\displaystyle\sqrt{r^{\mathsf T} r}$
		; \ req$_{r}$ := IAllReduce($r_{i}^{\mathsf T} B_{i}^{} r_{i}^{}$, $r^{\mathsf T} r$, SUM)
		; \ $k$ := $k + 1$
	\textbf{end if}
\ENDWHILE
\RETURN $x^{}_{i}$
\end{algorithmic}
}
\end{algorithm}
As one can see, compared to all the shared Boolean variables managed in the one-sided implementation (at both subdomain and coarse domain sides), non-blocking collective routines allow for a more straightforward implementation without any particular performance issue. In practice, the function \textit{IReduce} (applied to the contributions $R_{0}^{} R_{i}^{\mathsf T} B_{i}^{} r_{i}^{}$, $i \in \{1, \ldots, p\}$) can actually consist of an \textit{IGather} where the subdomains send only data related to their coarse unknown. It must be noted that, if an RDMA-capable network device is targeted, a hybrid approach is advisable, which consists of using one-sided routines for inter-subdomains communication and non-blocking collective routines for coarse-space correction.

Beyond the type of communication routines one chooses to consider, the main performance issue of such a two-level asynchronous solver is the weak incorporation of the coarse space. As suggested by \cite{GlusaEtAl2020b, GlusaEtAl2020}, each coarse solution is applied only once, however, if the coarse domain can just be seen as a $(p+1)$-th subdomain, then the corrections should just be done using the latest available coarse solution $x_{0}$, without wondering if it has been updated or not. In the sequel, we propose a more effective coarse-space correction where $x_{0}$ is more often incorporated.

\subsection{Accurate two-level asynchronous RAS solver}

The two-level asynchronous RAS solver \cite{GlusaEtAl2020b} features a one-sided non-blocking synchronization which consists of collecting all the residual contributions, $R_{0}^{} R_{i}^{\mathsf T} B_{i}^{} r_{i}^{}$, $i \in \{1, \ldots, p\}$, before computing the coarse residual, $r_{0}^{}$. Algorithm~\ref{algo:ras:async:cs} has provided a two-sided counterpart by using non-blocking collective routines. It however appears in both implementations that the computed coarse residual, $r_{0}^{}$, is inconsistent since each fine local component $r_{i}^{}$ is based on globally inconsistent interface data $x_{\neg i}^{}$. The effectiveness of the coarse-space correction can therefore be improved by considering a consistent coarse residual, which is achieved here by using consistent interface data $x_{\neg i}^{}$.

Let, then, $\bar x$ and $\bar r$ respectively denote the global solution and the global residual vectors that we shall consistently build. Algorithm~\ref{algo:ras:async:cs:isync} describes the improved two-level asynchronous solver.
\begin{algorithm}[htbp]
\caption{Accurate two-level asynchronous RAS.}
\label{algo:ras:async:cs:isync}
{\footnotesize
\begin{algorithmic}[1]
\REQUIRE $A_{i}^{}$, $R_{i}^{} A R_{\neg i}^{\mathsf T}$, $b^{}_{i}$, $x_{i}^{}$, $x_{\neg i}^{}$, $B_{i}^{}$, $N_{i}$, $A_{0}^{}$, $R_{0}^{} R_{i}^{\mathsf T}$, $i_{0}^{}$, $\varepsilon$, $\|b\|$, $k_{\text{max}}$
\STATE $r_{i}^{}$ := $b_{i}^{} - A_{i}^{} x_{i}^{} - R_{i}^{} A R_{\neg i}^{\mathsf T} x_{\neg i}^{}$
; \ AllReduce($r_{i}^{\mathsf T} B_{i}^{} r_{i}^{}$, $r^{\mathsf T} r$, SUM)
; \ $\|r\|$ := $\displaystyle\sqrt{r^{\mathsf T} r}$
; \ $k$ := 0
\STATE req$_{r}$ := REQUEST\_NULL
; \ state$_{0}$ := 0
; \ $x_{0}^{}$ := 0
\WHILE{$\|r\| \ge \varepsilon \|b\|$ \AND $k < k_{\text{max}}$}
	\STATE \textbf{if} state$_{0}$ = 0 \textbf{then}
		$\bar x_{i}^{}$ := $x_{i}^{}$
		; \ reqs$_{\bar x}$[ ] := ISynchronize($\bar x_{i}^{}$, $\bar x_{\neg i}^{}$, $N_{i}$)
		; \ state$_{0}$ := 1
	\textbf{end if}
	\IF{state$_{0}$ = 1 \AND TestAll(reqs$_{\bar x}$)}
		\STATE $\bar r_{i}^{}$ := $b_{i}^{} - A_{i}^{} \bar x_{i}^{} - R_{i}^{} A R_{\neg i}^{\mathsf T} \bar x_{\neg i}^{}$
		; \ req$_{r_{0}^{}}$ := IReduce($R_{0}^{} R_{i}^{\mathsf T} B_{i}^{} \bar r_{i}^{}$, $r_{0}^{}$, $i_{0}^{}$, SUM)
		; \ state$_{0}$ := 2
	\ENDIF
	\IF{state$_{0}$ = 2 \AND Test(req$_{r_{0}^{}}$)}
		\STATE \textbf{if} $i$ = $i_{0}$ \textbf{then}
			$x_{0}$ := Solve($A_{0}^{}$, $r^{}_{0}$)
		\textbf{end if}
		\STATE req$_{x_{0}^{}}$ := IBcast($x_{0}^{}$, $i_{0}^{}$)
		; \ state$_{0}$ := 3
	\ENDIF
	\STATE \textbf{if} state$_{0}$ = 3 \AND Test(req$_{x_{0}^{}}$) \textbf{then}
		state$_{0}$ := 0
	\textbf{end if}
	\STATE $x^{}_{i}$ := $x^{}_{i}$ + 0.5 $\times$ (Solve($A_{i}^{}$, $r^{}_{i}$) + $R_{i}^{} R_{0}^{\mathsf T} x_{0}^{}$)
	\STATE reqs$_{x}$[ ] := ISynchronize($x_{i}^{}$, $x_{\neg i}^{}$, $N_{i}$)
	; \ FreeAll(reqs$_{x}$)
	; \ $r_{i}^{}$ := $b_{i}^{} - A_{i}^{} x_{i}^{} - R_{i}^{} A R_{\neg i}^{\mathsf T} x_{\neg i}^{}$
	\STATE \textbf{if} Test(req$_{r}$) \textbf{then}
		$\|r\|$ := $\displaystyle\sqrt{r^{\mathsf T} r}$
		; \ req$_{r}$ := IAllReduce($r_{i}^{\mathsf T} B_{i}^{} r_{i}^{}$, $r^{\mathsf T} r$, SUM)
		; \ $k$ := $k + 1$
	\textbf{end if}
\ENDWHILE
\RETURN $x^{}_{i}$
\end{algorithmic}
}
\end{algorithm}
A step is added to Algorithm~\ref{algo:ras:async:cs}, prior to the coarse residual computation. Each subdomain $i \in \{1, \ldots, p\}$ records a version, $\bar x_{i}^{}$, of its local solution $x_{i}^{}$, and triggers synchronization with its neighbors on this particular version. The returned requests are then checked at each subsequent iteration. Once they are all completed, a globally consistent local residual, $\bar r_{i}^{}$, can be computed using the recorded local solution, $\bar x_{i}^{}$, and the corresponding synchronized interface solution, $\bar x_{\neg i}^{}$. The remaining non-blocking coarse-space correction procedure is then based on such consistent local residuals.

With such a more consistent coarse solution, $x_{0}^{}$, we also now let the coarse-space correction occur at several successive iterations, simply using the latest available $x_{0}^{}$. As we will see from numerical experiments, this can result in a significantly more effective two-level asynchronous solver.

\section{Experimental results}
\label{sec:exp}

\subsection{Problem and general settings}

Similarly to \cite{GlusaEtAl2020}, the Poisson's equation,
$
- \Delta u = g,
$
is considered for evaluating the practical performance of the proposed asynchronous additive coarse-space correction. The three-dimensional domain $[0, 1]^{3}$ is considered with a uniform source $g = 4590$ and $u = 0$ on the whole boundary. The domain is decomposed on each dimension, and the local matrices and vectors, $A_{i}^{}$, $R_{i}^{} A R_{\neg i}^{\mathsf T}$ and $b_{i}^{}$, $i \in \{1, \ldots, p\}$, are generated in each subdomain by P1 finite-element approximation. Both $A_{i}^{}$ and $A_{0}^{}$ are then Cholesky-factorized.

The compute cluster is homogeneously composed of nodes consisting of 175 GB of memory and two 20-cores processors at 2.1 GHz (40 cores per node), connected through an Omni-Path Architecture (OPA) network at 100 Gbit/s. Each core runs exactly one MPI process which also handles exactly one subdomain.

In the following, we show performances of the RAS solver \eqref{eq:ras} (referred to as ``1L-Sync''), its two-sided asynchronous execution from Algorithm~\ref{algo:ras:async} (``1L-Async''), the two-level RAS solver \eqref{eq:ras:cs} (``2L-Sync''), and its two-sided asynchronous execution from Algorithm~\ref{algo:ras:async:cs:isync} based on the proposed consistent coarse residual. We first considered applying each coarse solution only once (``1corr-2L-Async''), based on \cite{GlusaEtAl2020b, GlusaEtAl2020}, then we allowed for a maximum of five corrections using the same coarse solution (``5corr-2L-Async''). Overall execution times are compared for a stopping criterion $\|b - A x\| < \varepsilon \|b\|$ with $\varepsilon = 10^{-6}$. A final check is synchronously carried out after the iterations loop.

\subsection{Scaling}

We consider a weak scaling experiment, where both the size of the problem and the number of subdomains are increased by the same ratio. The size of the local problems is kept constant at approximately 20000 unknowns, while the size of the global coarse problem equals the number of subdomains. One of the involved processes handles the whole coarse domain, which, up to $p = 1600$ subdomains, does not yet raise a significant overhead issue, compared to the size of the subdomain. A discrete overlap of two mesh steps is considered throughout the experiment.

\begin{figure}[htbp]
\begin{center}
\includegraphics[scale=0.38]{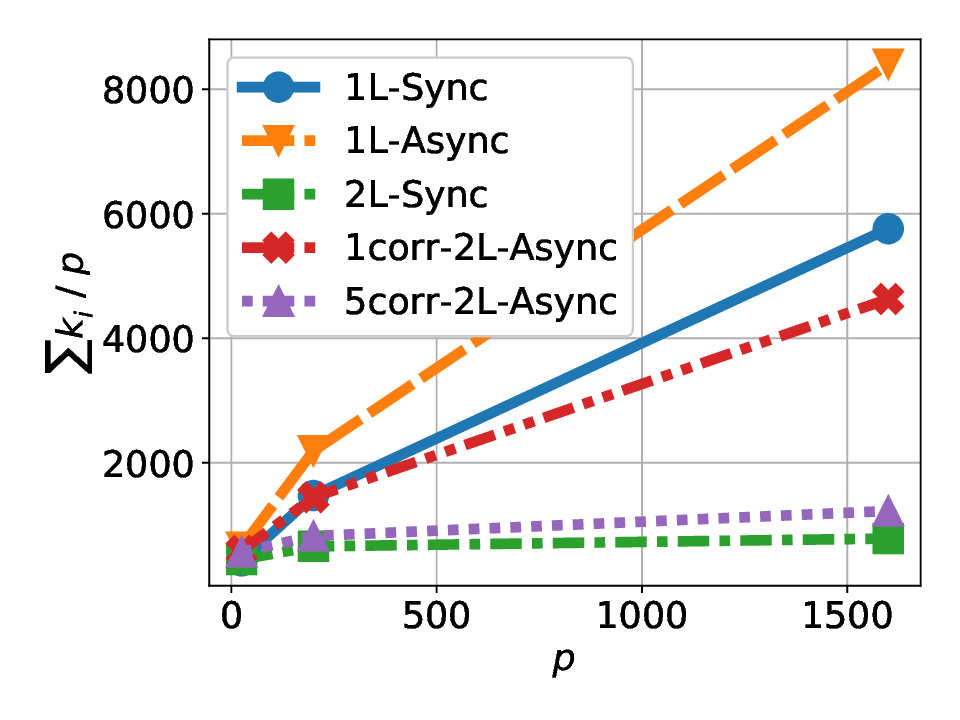}
\quad
\includegraphics[scale=0.38]{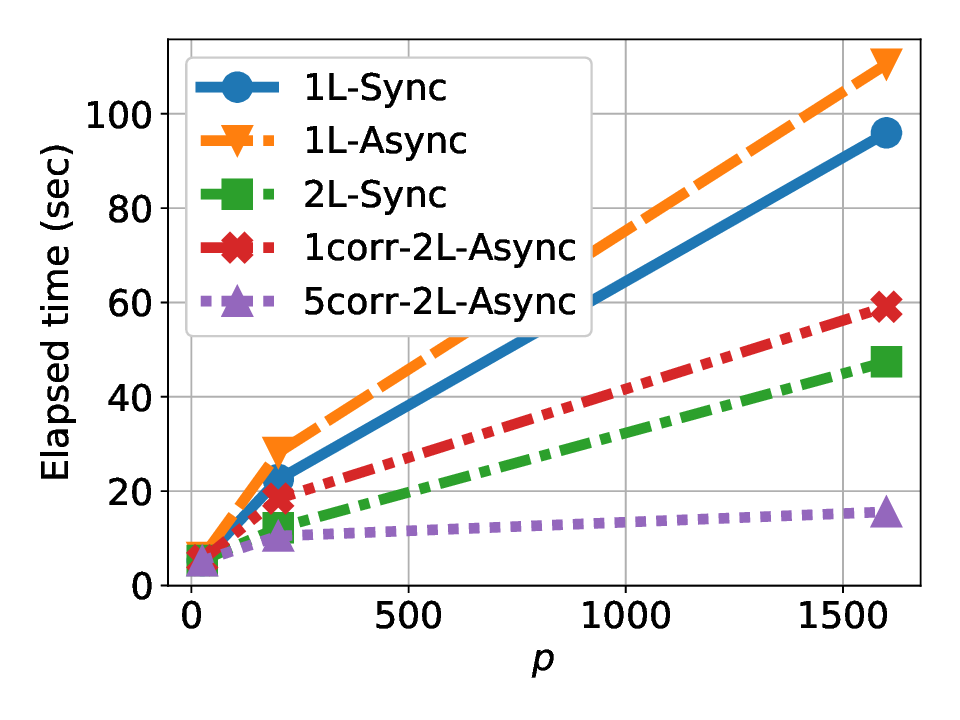}
\end{center}
\caption{Weak scaling of one-level and two-level RAS solvers. ``1corr-2L'' corresponds to 1 correction per coarse solution as in \cite{GlusaEtAl2020}, and ``5corr-2L'' allows for up to 5 corrections. $k_{i}^{}$ is the number of iterations on the $i$-th subdomain.}
\label{fig:scaling:weak}
\end{figure}

Numbers of iterations and total elapsed times are reported in Figure \ref{fig:scaling:weak}.
Comparing the synchronous solvers, there is a clear impact of the coarse-space correction. The number of iterations of the 2L-Sync solver is nearly constant, which confirms the effectiveness of the implemented coarse space. Comparison between synchronous and asynchronous solvers is relevant on elapsed time instead. For the one-level solvers, we note that, while showing the same scaling behavior, the synchronous one is still slightly faster at 1600 cores, which can be explained by the homogeneity of the cluster and the balanced workload. We observe the same result with the 2L-Sync and the 1corr-2L-Async solvers. Such a context, which is favorable to synchronous execution, strengthen the far better performance of the 5corr-2L-Async solver over the 2L-Sync and the 1corr-2L-Async solvers. This clearly confirms the benefit from allowing for more successive corrections using the same coarse solution, which is made possible here by improving the accuracy of the coarse residual. We should mention that the upper bound of 5 corrections per coarse solution became operative only at 1600 subdomains, where coarse solutions were less frequent, probably due to the cost of global communication towards the coarse domain.
This therefore strongly suggests the need for bounded delays in the theoretical analysis of two-level asynchronous iterations.

\section{Conclusion}
\label{sec:conclusion}

The development of effective asynchronous iterative solvers within the framework of domain decomposition methods is still at an early stage, however, promising results have been recently achieved, such as two-level asynchronous domain decomposition methods, within both overlapping Schwarz and primal Schur frameworks.
In this paper, we have proposed an improvement of the asynchronous additive coarse-space correction scheme applied in the overlapping Schwarz framework, which has resulted in a far more effective two-level asynchronous solver. Still, numerical experiments have also pointed out the necessity of deeper theoretical investigation of asynchronous coarse-space correction models, notably relying on bounded delays.

It turns out that effective implementation of asynchronous coarse space is an issue similar to the asynchronous convergence detection one, both of them being related to the accurate evaluation of a global residual. While accurate convergence detection is sometimes seen as of a low importance, building an accurate global residual now appears as a key ingredient for achieving effective two-level asynchronous solvers.

\section*{Acknowledgements}

This work was partly funded by the French National Research Agency as part of project ADOM, under grant number ANR-18-CE46-0008.

\bibliography{ref}
\bibliographystyle{abbrv}

\end{document}